\documentclass[%
 aip,
cp,  
 amsmath,amssymb,
 reprint,%
]{revtex4-2}

\usepackage{graphicx}
\usepackage{dcolumn}
\usepackage{bm}

\usepackage[utf8]{inputenc}
\usepackage[T1]{fontenc}
\usepackage{mathptmx} 
\renewcommand{\thetable}{\arabic{table}}    

\usepackage{amsmath}
\usepackage{mathtools}

\DeclareMathOperator*{\argmin}{arg\,min}

\newcommand{\bx}{{\bf x}}

\newcommand{\bA}{{\bf A}}
\newcommand{\bB}{{\bf B}}
\newcommand{\bC}{{\bf C}}

\newcommand{\bI}{{\bf I}}
\newcommand{\bK}{{\bf K}}
\newcommand{\bL}{{\bf L}}

\newcommand{\bQ}{{\bf Q}}

\newcommand{\bAw}{{\bA_\text{w}}}

\newcommand{\bKw}{{\bK_\text{w}}}
\newcommand{\bKL}{{\bK_\text{L}}}
\newcommand{\bAc}{{\bA_\text{c}}}
\newcommand{\bQc}{{\bQ_\text{c}}}
\newcommand{\bKc}{{\bK_\text{c}}}
\newcommand{\Cobs}{{C_\text{obs}}}
\newcommand{\Kobs}{{K_\text{obs}}}

\newcommand{\Kw}{{K_\text{w}}}
\newcommand{\Kc}{{K_\text{c}}}

\usepackage{bm} 
\newcommand{\bmxi}{{\bm{\xi}}}
\newcommand{\bmeta}{{\bm{\eta}}}
\newcommand{\bmtau}{{\bm{\tau}}}

\begin{document}

\title{A Liang-Kleeman Causality Analysis based on Linear Inverse Modeling}

\author{Justin Lien}  
\email[ Corresponding author:  ]{lien.justin.t8@dc.tohoku.ac.jp}  
\affiliation{ 
  Mathematical Institute, Tohoku University. \\
  6-3 Aramaki Aza, Aoba, Aoba Ward, Sendai, Miyagi, 980-0845, Japan.
}



\date{\today} 

\begin{abstract}
Causality analysis is a powerful tool for determining cause-and-effect relationships between variables in a system by quantifying the influence of one variable on another. 
Despite significant advancements in the field, many existing studies are constrained by their focus on unidirectional causality or Gaussian external forcing, limiting their applicability to complex real-world problems. 
This study proposes a novel data-driven approach to causality analysis for complex stochastic differential systems, integrating the concepts of Liang-Kleeman information flow and linear inverse modeling.
Our method models environmental noise as either memoryless Gaussian white noise or memory-retaining Ornstein-Uhlenbeck colored noise, and allows for self and mutual causality, providing a more realistic representation and interpretation of the underlying system. 
Moreover, this LIM-based approach can identify the individual contribution of dynamics and correlation to causality. 
We apply this approach to re-examine the causal relationships between the El Ni\~{n}o-Southern Oscillation (ENSO) and the Indian Ocean Dipole (IOD), two major climate phenomena that significantly influence global climate patterns.
In general, regardless of the type of noise used, the causality between ENSO and IOD is mutual but asymmetric, with the causality map reflecting an ENSO-like pattern consistent with previous studies. 
Notably, in the case of colored noise, the noise memory map reveals a hotspot in the Niño 3 region, which is further related to the information flow.
This suggests that our approach offers a more comprehensive framework and provides deeper insights into the causal inference of global climate systems. 
\end{abstract}

\maketitle

\section{Introduction}

Understanding causal relationships is fundamental to controlling the behavior of complex dynamical systems and developing effective solutions to intricate real-world problems \cite{Guo2020}. 
Rigorous causality analysis allows us to move beyond mere descriptions of phenomena towards a deeper, explanatory understanding of the underlying processes driving observed outcomes. 
Applications of causality analysis span diverse fields, from climate science to machine learning, enabling the optimization of business strategies and technological systems \cite{Alper2016,Ye2022}. 
Identifying and quantifying causal effects is indispensable for advancing scientific knowledge and driving impactful innovation.

While substantial progress has been made in causal discovery, many existing techniques struggle to handle complex systems. 
A prominent example is the Direct LiNGAM, which effectively identifies causal structures in \textit{static systems} subjected to identically and independently distributed non-Gaussian noise with unidirectional causality \cite{Majdi2024,Shimizu2011}.
However, real-world phenomena are often characterized by complex dynamics where feedback loops and bidirectional or mutual causal interactions among variables are ubiquitous. 
Such temporal characteristics and intricate causal dependencies cannot be fully captured by the Direct LiNGAM framework.

To address these limitations, several studies have explored Liang-Kleeman information flow theory, which provides a quantitative measure of directed information transfer and offers a more dynamic perspective on causality \cite{Liang2005,Liang2007,Liang2008}. 
Building upon this foundation, Liang proposed an empirical method that identifies causal relations among variables sampled from a stochastic dynamical system \cite{Liang2014,Liang2021}. 
Contrary to the "static one-way causality" assumed by Direct LiNGAM, Liang's approach includes the presence of self and mutual causality, and captures the dynamical interactions between variables, where some variables may stabilize or excite others. 
However, it models the random forcing with Gaussian white noise, which lacks memory and temporal dependency. While this assumption simplifies the mathematical modeling, it may not accurately represent many real-world systems where the random forcing exhibits temporal correlations.

Colored noise, a type of temporally correlated noise with an inherent memory, has gained attention in mathematical modeling as it can provide a more realistic representation of environmental noise \cite{Jin2007,Zhao2024}.
Linear Inverse Models (LIMs), a class of data-driven approaches that extract dynamical and random forcing information from observations, have been extended to allow for Gaussian white and Ornstein-Uhlenbeck (OU) colored noise \cite{Lien2024,Penland1989}. 
It has been shown that modeling environmental random forcing with colored noise can alter both deterministic dynamics and stochastic forcing, suggesting that the memory effect of the noise can also have significant implications for causal inference \cite{Jin2007}.

In this paper, we propose a new method for causality analysis that accommodates both white and colored noise by combining the concepts of information flow and LIMs. 
Our proposed method can determine causality in dynamic systems, accounting for self and mutual interactions and the memory effect of the environmental noise. 
In addition, the individual contribution from dynamics and correlation to causality can be identified and quantified. 
This generalization provides a more comprehensive and realistic tool for causality analysis in a wide range of dynamical systems.

\section{Methodology}
We start with the general idea. Let's suppose that the $n$-dimensional observation data is collected from 
\begin{align} \label{Eq:GE} 
    \frac{d}{dt} \bx = F(\bx,t,\text{noise})
\end{align}
where $F$ describes the complex system of interest.
Generally speaking, while the external stochastic noise continuously excites the state variables $\bx_i$'s, the deterministic dynamics act as a damping operator, dissipating energy and drawing the variables back to equilibrium \cite{Penland1989}.
We aim to identify the causal inference between $\bx_i$'s, but
the system $F$ is mathematically intractable in practice.
Therefore, to extract the information of the underlying unknown system, an equation-free perspective is taken by constructing a linear system to approximate Eq. (\ref{Eq:GE}):
\begin{align} \label{Eq:AE} 
    \frac{d}{dt} \bx = \bA \bx + \sqrt{2\bQ} \cdot \text{noise} 
\end{align}
where $\bA$ is the dynamical matrix whose eigenvalues have negative real parts, and $\bQ$ is the positively definite diffusion matrix.
Then the problem of identifying the causality is formulated as an inverse problem: given an observation time series, find a suitable linear representation (\ref{Eq:AE}) whose causality estimates the causality of the underlying unknown system (\ref{Eq:GE}).

Previous studies in the inverse problem mostly model the environmental forcing as the memoryless normalized white noise $\bmxi$ while, in this study, the noise term is allowed to be OU colored noise $\bmeta$ \cite{Penland1989,Lien2024}. 
This is the steady state solution of the OU process $\frac{d}{dt}\bmeta = \frac{-1}{\bmtau}\bmeta + \frac{1}{\bmtau}\bmxi$ whose temporal correlation is of an exponential shape: $\langle \bmeta(t) \bmeta(s)^T \rangle = \frac{1}{2\bmtau}e^{\frac{-\left\vert t-s \right\vert}{\bmtau}}\bI$, where $\bI$ denotes the identity matrix and the bracket denotes the expectation. 
The parameter $\bmtau$, called the noise correlation time, can be physically interpreted as the $e$-folding time of the noise temporal dependency.  
To avoid confusion whenever necessary, we use the script $\cdot_\text{w}$ and $\cdot_\text{c}$ to distinguish the choice of the noise term.
In the rest of this section, we focus on the definition of causality and the construction of the linear representation (\ref{Eq:AE}) of the unknown system (\ref{Eq:GE}).

\subsection{Liang-Kleeman information flows}
The rate of information flowing from a component $\bx_j$ to another $\bx_i$, denoted by $T_{j\to i}$ (unit: nats/time), is defined to be the change rate of the marginal (Shannon) entropy of $\bx_i$, minus the same change rate but without the effect from $\bx_j$ to the system \cite{Liang2008}.
With the information flow, the Liang-Kleeman causality is quantitatively defined: (1) $\bx_j$ is causal to $\bx_i$ if and only if $T_{j\to i} \ne 0$, and (2) the magnitude of the causality from $\bx_j$ to $\bx_i$ is measured by $\left\vert T_{j\to i} \right\vert$ \cite{Liang2021}.
In general, the information flow admits a closed form for any stochastical differential equation (SDE) with appropriate regularity but it is considerably complicated and out of the scope of this study \cite{Liang2014}.

For the linear system with white noise random forcing, Liang shows that the information flows $T_{j \to i}$ from $\bx_j$ to $\bx_i$ is given by
\begin{align} \label{Eq:IF0}
    T_{j \to i} = \bA_{ij} \frac{\bC_{ij}}{\bC_{ii}},
\end{align}
where $\bC = \langle \bx(\cdot) \bx(\cdot)^T \rangle$ is the covariance of the state variables \cite{Liang2021}.
Now the physical interpretation of information flow and hence causality becomes clear. 
Without loss of generality, we may assume that $\bx_i > 0$ at some moment.  
If $\bC_{ij} > 0$ and $\bA_{ij} > 0$, then the $\bx_j$ is more likely to have the same sign as $\bx_i$ and the contribution from $\bx_j$ to $\bx_i$ turns out to be positive, driving the $\bx_i$ away from the steady state and making it unstable.
By considering all possible combinations of signs, we can easily see that $T_{j \to i} > 0$ ($<0$) means that $\bx_j$ excites (stabilizes) $\bx_i$.
We note that Eq.~(\ref{Eq:IF0}) is available for self-causality $T_{i \to i}$, which solely depends on the dynamics $\bA_{ii}$.

On the other hand, for the linear system with colored noise random forcing, the governing equation can be written as the white-noise-driven augmented system
\begin{align} 
    \frac{d}{dt}
    \begin{bmatrix}
        \bx \\
        \bmeta
    \end{bmatrix}
    = 
    \underbrace{
    \begin{bmatrix}
        \bAc &  \sqrt{2\bQc} \\
        0 & -\frac{1}{\bmtau}\bI 
    \end{bmatrix}
    }_{\bL}
    \begin{bmatrix}
        \bx \\
        \bmeta
    \end{bmatrix}
    +
    \frac{1}{\bmtau}
    \begin{bmatrix}
        0 \\
        \bmxi 
    \end{bmatrix}
    ,
\end{align}
where $\bL$ plays the role of augmented dynamics.
We notice that the deterministic interaction among state variables $\bx_i$'s remains unchanged. 
Therefore, Eq. (\ref{Eq:IF0}) is still valid for the colored-noise case, and so is the physical interpretation of the information flow.

When the noise term is white, Liang used the maximum likelihood estimation to show that Eq. (\ref{Eq:IF0}) can be expressed as a function of the covariance function of the state variables and their derivative as follows,
\begin{align} \label{Eq:IF1}
    T_{j \to i} = \frac{1}{\det \bC} \sum_{k = 1}^n \Delta_{jk} \bC_{k,di} \cdot \frac{\bC_{ij}}{\bC_{ii}}
\end{align}
where $\Delta_{jk}$ is the cofactor of $\bC$ and $\bC_{k,di}$ is the covariance between $\bx_k$ and $\Dot{\bx_i}$ \cite{Liang2021}.
Therefore, given observation data, the information flows $T_{j \to i}$ can be estimated by Eq. (\ref{Eq:IF1}) in which $\Dot{\bx_i}$ can be computed by forward finite difference scheme in practice.
However, using Eq. (\ref{Eq:IF1}) causes two major limitations.
Firstly, the correlation-relevant information among state variables is mixed, and it becomes difficult to identify their individual roles in the causality. 
Moreover, in the case of colored noise, the covariance matrices of $\bx$ (and its derivative) in Eq. (\ref{Eq:IF1}) should be replaced by those of $[\bx;\bmeta]$, the concatenation of the state variables and colored noise.
Therefore, Eq. (\ref{Eq:IF1}) is infeasible for the inverse problem since the colored noise is not observed in practice.
To overcome these limitations, we employ the Colored-LIM to estimate the linear dynamics $\bAc$ without requiring direct observation of $\bmeta$ so that a direct application of Eq. (\ref{Eq:IF0}) is feasible.

\subsection{Linear inverse modeling}

The linear inverse model is a class of data-driven methods that extracts the dynamic and stochastic relevant information from observation data. 
More precisely, in the LIM framework, the unknown complex system (\ref{Eq:GE}) is approximated by the linear system (\ref{Eq:AE}). 
The dynamical matrix $\bA$ of such a linear system is encoded in the correlation function of the state variables, and the diffusion matrix $\bQ$, though not necessary for this study, follows from the fluctuation-dissipation relations \cite{Penland1989,Lien2024}. 
In this study, we focus on White-LIM and Colored-LIM, in which the random forcing is modeled by Gaussian white noise and OU colored noise, respectively.

\subsubsection{White-LIM}
It is well-known that the white-noise-driven linear system is an OU process whose correlation function $\bKw = \langle \bx(\cdot + s)\bx(\cdot)^T\rangle$ is an exponential function,
\begin{align} \label{Eq:bKw}
    \bKw(s) = e^{s\bAw} \bC
\end{align}
for $s \ge 0$ \cite{Penland1989}. 
The exponent of $\bKw$ is the linear dynamics that can be obtained by using the matrix logarithm
\begin{align} \label{Eq:Aw}
    \bA_\text{w}(\rho) = \frac{\log(\bKw(\rho)\bC^{-1})}{\rho}
\end{align}
for any arbitrary time-lag $\rho$.
In practice, due to the sampling size and model uncertainty (e.g., non-linearity, non-Gaussianity), the estimated dynamics $\bAw$ may be unstable (i.e., $\bAw$ strongly depend on $\rho$), and the choice of $\rho$ may differ from studies \cite{Penland1995}.
Indeed, we observe that Eq. (\ref{Eq:Aw}) is equivalent to finding an exponential function that passes through the correlation function at $0$-lag and $\rho$-lag.
Hence, the instability may be due to the insufficient data used to determine $\bAw$.

\subsubsection{Colored-LIM}
For a colored-noise-driven linear system, the correlation function $\bKc$ reads, for $s \ge 0$,
\begin{align} \label{Eq:Kc}
    \bKc(s) = e^{s\bAc}\bC + e^{s \bAc} \int_0^s e^{-s'(\bAc+\frac{1}{\bmtau})} ds' \, \bQc\bB^T
\end{align}
where $\bB=(1-\bmtau\bAc)^{-1}$ accounting for the noise memory. 
The first term is the exponential decay due to the deterministic dynamics while the second term involves the contribution from the state-dependent colored noise.
As Eq. (\ref{Eq:Kc}) is no longer an exponential function, Eq. (\ref{Eq:Aw}) is infeasible. 
Nevertheless, Lien et al. show that with $\bmtau$ treated as a hyperparameter, the higher derivatives of $\bKc$ encode the linear dynamics $\bAc$ which can be obtained by solving linear equations \cite{Lien2024}. 
This simplifies the problem and speeds up the computation time but a robust algorithm for numerical differentiation is required.
Though in the real-world application, colored noise may be a better representation of the environmental noise, the derivative-based approach sacrifices the degree of freedom for $\bmtau$.

\subsection{LIM-based algorithm for causal analysis}
In this study, we propose a LIM-based approach for Liang-Kleeman causality analysis by utilizing the ideas of White-LIM and Colored-LIM to obtain the linear dynamics $\bA$ and then applying Eq. (\ref{Eq:IF0}) to estimate the information flow. 
Here, we merely focus on the estimation of $\bA$.

Given observation data, let $\Cobs$ and $\Kobs$ denote the observed covariance and correlation function, respectively. To avoid the numerical instability in the White-LIM and the treatment of $\bmtau$ as a hyperparameter in the Colored-LIM, we resort to a minimization method to compute linear dynamics:
\begin{align} 
    \bAw = \argmin_{A} \left\Vert (\Kobs - \Kw)|_{[0,l]} \right\Vert_{w, F} \label{Eq:3a} 
    \intertext{for the white-noise case and}
    (\bAc, \bmtau) = \argmin_{A,\tau} \left\Vert (\Kobs - \Kc)|_{[0,l]} \right\Vert_{w, F} \label{Eq:3b}
\end{align}
for the colored-noise case, under the constraint that $\Kw(0) = \Kc(0) = C_\text{obs}$, where $l$ is the window length that specifies the minimization window, $\left\Vert \cdot \right\Vert_{w, F}$ denotes the weighted Frobenius norm, and $\Kw$ and $\Kc$ are the correlation functions determined by the running variables $A$ (and $\tau$), Eqs.~(\ref{Eq:bKw}) and (\ref{Eq:Kc}).
In the White-LIM framework, the minimization process takes several $\Kobs(s)$ with $s > 0$ into account instead of simply choosing one of them, while in the Colored-LIM framework, the noise correlation time is also a variable in the search space, opening the possibilities to study the noise memory.
    
The minimization process requires a proper initial guess for a robust result. 
Though Eq. (\ref{Eq:Aw}) may not be stable, it returns the linear dynamics $A_0(\rho)$ which determines a correlation function that at least passes two points of $\Kobs$. 
As a result, $A_0(\rho)$ for different $\rho$'s are educated initial guesses.
As for the initial guess for $\tau$ in the Colored-LIM framework, since in the limit $\tau \to 0$, the colored noise formally approaches the white noise, we can start from a small $\tau$ and then adjust it to optimize the minimization process.

We remark that the White-LIM and Colored-LIM-based approaches explicitly quantify the contribution of dynamics, which cannot be directly inferred by using Eq. (\ref{Eq:IF1}).
Moreover, though the LIM-based approaches
use the same formula (\ref{Eq:IF0}) to estimate the information flow and hence the causal implication, the linear dynamics indeed depends on the choice of noise, especially when the environmental noise is non-trivially temporally correlated.
As a result, they can lead to different information flows and causality in both sign and magnitude, even if applied to the same observation data.

\section{A real-world application}

We apply the LIM-based algorithm for causal analysis to study the causal relation between two major climate phenomena the El Niño-Southern Oscillation (ENSO) and the Indian Ocean Dipole (IOD). 
ENSO is an interannual climate phenomenon involved in the ocean-atmosphere interaction over the tropical Pacific Ocean \cite{Cai2015}. 
It has been shown that ENSO can influence local weather patterns, increasing disaster risks like droughts and floods in an extreme ENSO event \cite{Cai2015,Philip2014}.
Detecting extreme ENSO events requires observations from both the atmosphere and oceans, but the average sea surface temperature (SST) over several regions of the tropical Pacific Ocean has been identified as being important for monitoring such events. 
Examples include the Ni\~{n}o 3 region (150W-90W and from 5S-5N). 
The IOD, on the other hand, is an irregular large-scale oscillation of SST and atmospheric circulation in the Idena Ocean \cite{Saji1999}. 
Similar to ENSO, an extreme IOD event is linked to extreme weather events in Indian Ocean rim countries \cite{MacLeod2024}.
IOD is commonly measured by the Dipole Mode Index (DMI), the SST difference between the western (50E-70E and 10S-10N) and the south-eastern equatorial Indian Ocean (90E-110E and 10S-0N).
These two systems can interact with each other, but the teleconnection is still unclear \cite{Hrudya2021}.
Understanding their interrelationship helps predict and prepare for significant weather changes affecting agriculture, water resources, and economies globally \cite{Reddy2022}.

In the LIM framework, the state variables are representative of the studied climate phenomenon, where the deterministic dynamics is often regarded as the oceanic systems or large-scale atmosphere-ocean systems, but the stochastic noise may have a different physical meaning \cite{Matthew2011}.
White noise, which is widely applied in most LIM applications, is considered as atmospheric random events such as weather disturbances and precipitation \cite{Penland1995}. 
Colored noise, on the other hand, can be thought of as modeling the related climate or weather events such as Madden-Julian Oscillation and westerly wind bursts, or the integrated effect of climate variables other than the state variables \cite{Zhao2024}.
In this study, we collect the DMI and Pacific SSTs data from NOAA Physical Sciences Laboratory \cite{NOAA2024,Huang2017} to be the state variables.
To compare our results with those of Liang's study, we use data from the same period (from 1958 to 2010) for the following discussion.
Following the convention in the climate science community, we remove the seasonal variation of observation data and then apply a $3$-month running average. 
The minimization window length $l$ is set to be $4$ months.

\begin{figure}
\includegraphics[]{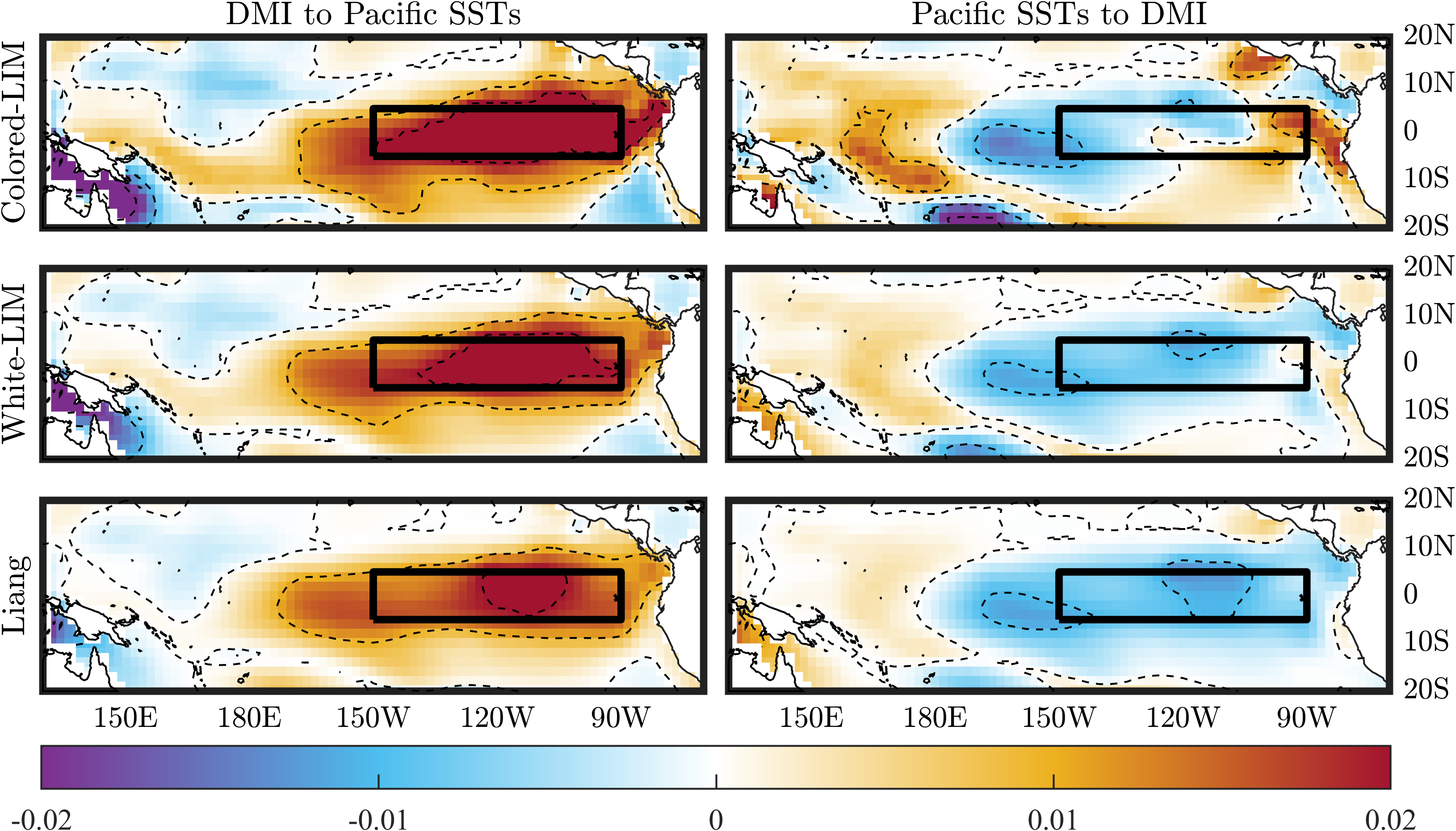}
\caption{\label{Fig:1} The distribution of information flows (unit: nats/months) between DMI and Pacific SSTs computed by LIM-based methods, and Liang's method. The information flows larger than $0.02$ (smaller than $-0.02$) are masked out to $0.02$ (-0.02). The block box specifies the Ni\~{n}o 3 region. }
\end{figure}

Figure \ref{Fig:1} shows the information flows between DMI and each grid point in the mid-Pacific ocean. 
Though slightly different, the distribution of information flows for each method shows an El Ni\~{n}o-like pattern.
It is clear that the causality is mutual but asymmetric: Pacific SST tends to stabilize IOD while IDO excites Pacific SST.
The similarity between Liang's method and our White-LIM-based approach is not surprising, as both approximate the underlying system as a white-noise-driven linear process. 
However, we notice that the Colored-LIM-based results tend to reveal a stronger information flow, especially over the Ni\~{n}o 3 region. 
In general, for fixed observation data, when the noise memory $\bmtau$ is non-trivial, the external random forcing $\bmeta$ is more likely to stay in the same sign over a period of time, consistently driving the state variables $\bx_j$'s to the extreme. 
However, as the climate system mitigates the occurrence of extreme events (corresponding to the fixed $\Cobs$ in the LIM framework), the damping force $\bAc$ must be stronger to fight against the effect of memory noise.
Consequently, the entry of $\bAc$ turns out to be more likely to be larger in magnitude, and so is the information flow. 

\begin{figure}
\includegraphics[]{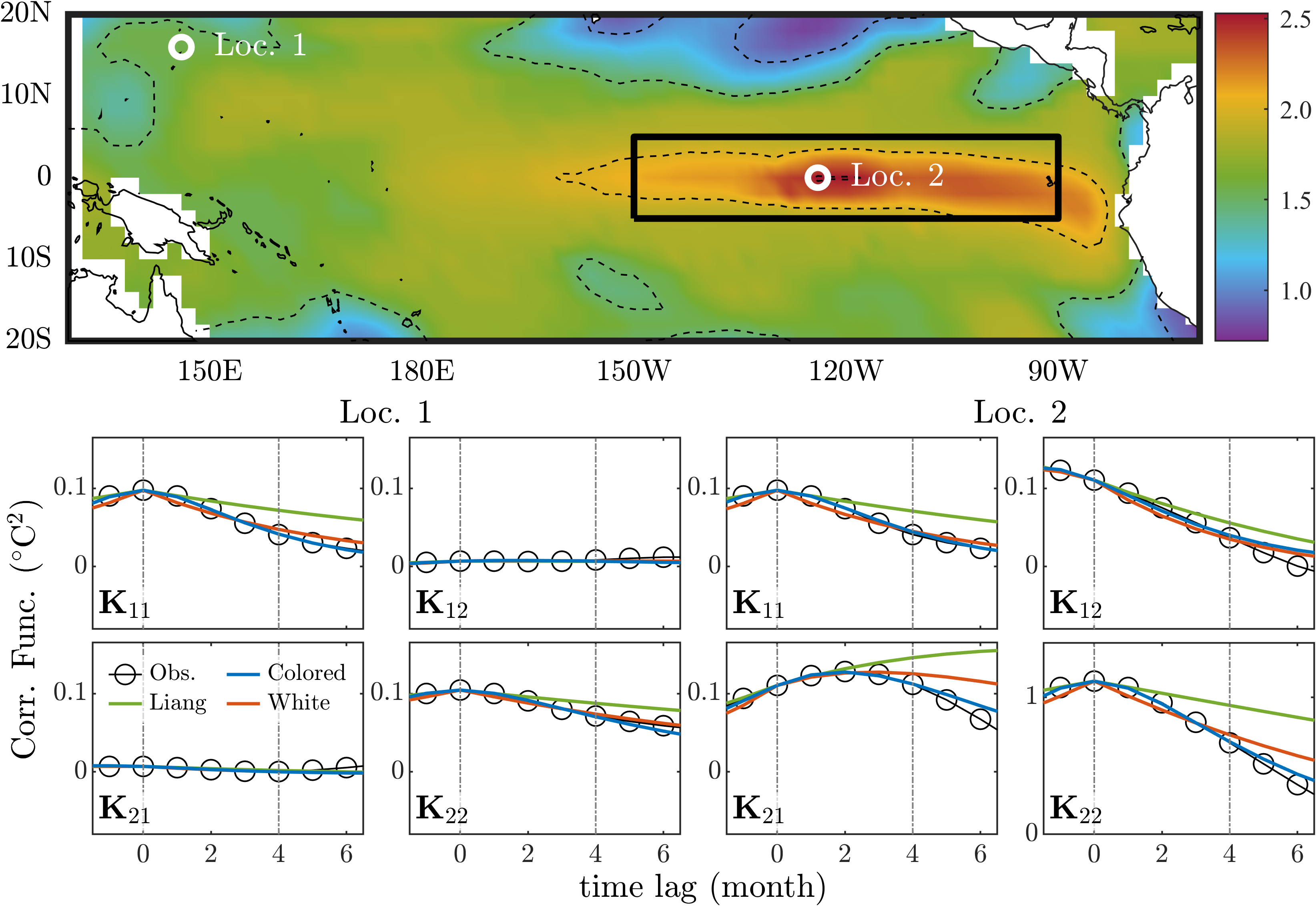}
\caption{\label{Fig:3} The distribution of $\bmtau$ (upper panel, unit in months), and $(i,j)$-entry of the observed and theoretical correlation functions (lower panels). The white circles in the upper panel specify the locations from which data is used to compute the correlation functions in the lower panels. The DMI and Pacific SSTs are symbolized by scripts 1, and 2, respectively. The dashed line represents the window over which the minimization is applied. 
The theoretical correlation function of Liang's method is computed by using Eqs. (\ref{Eq:IF0}), (\ref{Eq:IF1}), and (\ref{Eq:bKw}).
The scale of the $y$-axis for SST auto-correlation ($\bK_{22}$) in Loc. 2 is significantly larger than others. }
\end{figure}

Contrary to Liang's method and the White-LIM-based algorithm, the Colored-LIM-based algorithm further reveals the spatial distribution of noise memory $\bmtau$.
Figure \ref{Fig:3} shows that the noise memory over the Ni\~{n}o 3 region reaches more than 2 months, approximately 50\% longer than the rest of the mid-Pacific ocean. 
To understand the implication, we focus on the observed correlation functions $\Kobs$ and the theoretical correlation function $\bKw$ and $\bKc$ determined by Eqs.~(\ref{Eq:bKw}), (\ref{Eq:Kc}), (\ref{Eq:3a}), and (\ref{Eq:3b}) at two locations: Loc. 1 (146E and 16N) away from the Ni\~{n}o 3 region and Loc. 2 (124W and 0N) at which the maximum $\bmtau$ is achieved.

First, the observed auto-correlations for DMI and SSTs (i.e., the diagonal entries of $\Kobs$) at both locations are concave downward at the origin.
This is a characteristic of temporally correlated random forcing in the linearized framework, since the noise memory makes the trajectory of the state variables more sticky, making $\bK''(0)$ negatively definite \cite{Lien2024}. 
In contrast, the memoryless white noise causes a more oscillating trajectory and hence, a concave upward feature at the origin.
Moreover, at Loc. 1, the observed auto-correlation of SST is relatively small, and cross-correlation (i.e., the off-diagonal entries of $\Kobs$) is nearly constant and negligible, meaning that the state variables can be considered as decoupled, implying that the 1-month noise memory is internal for each ocean. 
At Loc. 2, the auto-correlation for SST is significantly larger, reflecting the strong climate variability of ENSO, and the time-dependent cross-correlation is more prominent than that at Loc. 1, suggesting a stronger correlation between these two climate phenomena. 
Within the Colored-LIM framework, whenever the noise memory is long, the damping force must be stronger, due to the continuous excitation of the colored noise as mentioned before.
Besides, when a stronger time-dependent cross-correlation between variables exists in the system, their interaction or coupling (i.e., the off-diagonal elements of $\bAc$) is more likely to be stronger in response to fit the observation.
As a result, the information flows over the Ni\~{n}o 3 region is comparably significant to other parts of the Mid-Pacific Ocean, suggesting a direct connection among noise memory, dynamics, cross-correlation, and causality.

Finally, we reconstruct the white-noise-driven linear system described in Eq. (\ref{Eq:AE}) implicitly used in Liang’s method, and show its correlation functions $\bKL$ in Figure \ref{Fig:3}. 
While it is evident that $\bKL$ aligns well with observations for small lags $s \le 2$ months, the deviation rapidly increases for larger lags. 
This occurs because Eq. (\ref{Eq:IF1}) only considers the covariance of neighboring observations (i.e., $0$-lag and $\Delta t$-lag, where $\Delta t$ denotes the sampling interval). In contrast, the theoretical correlations $\bKw$ and $\bKc$ derived from the LIM-based methods show strong alignment with the observed data, owing to the minimization technique employed. 
Although it remains challenging to identify the optimal linear system (\ref{Eq:AE}) that best represents the complex earth system, it is understandable that the three methods yield slightly different patterns for Pacific-SST-to-IOD causality. 
Nevertheless, these different models can offer valuable and distinct insights into the same problem.

\section{Concluding remarks}

In this study, we have proposed a LIM-based approach to causal analysis for complex dynamical systems. 
This method provides an additional option for modeling the external forcing by OU colored noise, and also allows the separation of the influence of dynamics and correlation on causality. 
Indeed, from observation, we can merely discuss the correlation between variables, but with the help of LIM, we can further discuss the dynamics. 
Finally, utilizing the information flow, we arrive at the causal inference.

Applying the LIM-based approach to the DMI and Pacific SSTs, we verified that they are mutually causal but asymmetric, and the distribution of information flow reveals an El Ni\~{n}o-like pattern, consistent with Liang's result. 
Besides, we further show that the environmental noise of the climate system at the Ni\~{n}o 3 region has a memory of $\bmtau \gtrsim 2$ months, longer than other parts of the Pacific Ocean. 
Such a noise memory further enhances the deterministic damping force and the coupling, leading to stronger information exchanged between the Ni\~{n}o 3 region and Indian Oceans.

Liang's method (\ref{Eq:IF1}) is closely related to our White-LIM-based method.
The right derivative of Eq. (\ref{Eq:bKw}) at the origin is $\langle \Dot{\bx}(\cdot) \bx(\cdot)^T \rangle = \bAw \bC$, and the first half on the right-hand side of Eq. (\ref{Eq:IF1}) is equivalent to the entry-wise expression of $\bAw = \langle \Dot{\bx}(\cdot) \bx(\cdot)^T \rangle \bC^{-1}$.
In practice, the computation of $\Dot{\bx}$ can be done by the forward difference scheme.
This is equivalent to extracting the linear dynamics by finding a straight line passing through $\Kobs$ at the origin and $\Delta t$-lag. 
In contrast, the White-LIM achieves the same task by finding an exponential function that passes through the origin and an arbitrary lag through Eq. (\ref{Eq:Aw}). 
These correspond to the linear and exponential fitting in \cite{Lien2024b}. 
Therefore, if the underlying system is a white-noise-driven linear system, the White-LIM-based approach is expected to reach a better performance by avoiding the numerical error in the finite difference scheme.

As in Liang's method where one is free to choose a finite difference scheme, the LIM-based approaches have a degree of freedom: the window width $l$ for minimization. 
Based on our experience, the numerical results are more dependent on $l$ when it is small due to the same reason as the numerical instability of the White-LIM. 
On the other hand, a larger $l$ does not necessarily improve the performance, since the correlation at large $l$ may be strongly influenced by, for example, the non-linearity of the underlying system, rendering the linear approximation (\ref{Eq:AE}) ineffective. 
Therefore, a mid-range $l$ should be appropriate.
We have also checked that $l = 5$ does not significantly change the discussion in this study. 

With the same idea, our LIM-based framework can also be applied to periodic systems. 
In fact, the cyclo-stationary LIM (CS-LIM) is a variant of LIMs that estimates the periodic dynamics by interval-wise approximating the underlying system with Eq. (\ref{Eq:AE}) from a cyclo-stationary time series \cite{Lien2024b,Shin2021}.
In the climate sciences, it is shown that such a method better represents the seasonal variation and enhances the prediction skills \cite{Shin2021}.
Utilizing the CS-LIM potentially enables us to study the seasonal dependency of information flows, leading to a more detailed temporal pattern of causality between oceans.

Finally, we remark that the study of the climate system has not yet been completed. 
For instance, the physical meaning of the colored noise and its memory is still to be clarified. 
The pathway of colored noise to the oceans requires further investigation.
Also, the positive Pacific-SST-to-IOD causality near South America in the Colored-LIM-based approach seems intriguing. 
These will be left for future study.

\begin{acknowledgments}
The author would like to express gratitude to NOAA PSL for making NOAA Extended Reconstructed SST V5 data publicly available \url{https://psl.noaa.gov}, and Professor Hiroyasu Ando for his constructive feedback and insightful discussions.
\end{acknowledgments}

\nocite{*}
 \bibliography{aipsamp} 

\end{document}